\documentclass[12pt]{article}

\usepackage{xcolor}
\usepackage[english]{babel}
\usepackage{amsmath,amsthm}
\usepackage{amsfonts}
\usepackage{float}
\usepackage{indentfirst}
\usepackage{subcaption}
\usepackage[hang,flushmargin]{footmisc}
\usepackage{amssymb}

\allowdisplaybreaks[4]

\usepackage{titling}

\pretitle{\begin{tabular}\large}
\posttitle{\par\end{tabular}}

\preauthor{\begin{center}\normalsize\begin{tabular}{c}}
\postauthor{\end{tabular}\par\end{center}}

\usepackage{graphicx}
\newtheorem{thm}{Theorem}
\newtheorem{claim}{Claim}
\newtheorem{subclaim}{Claim}[claim]
\newtheorem{lem}{Lemma}
\newtheorem{cor}{Corollary}

\newtheorem{proof of Thm}{proof of Thm}

\pretitle{\begin{center}\large}
\posttitle{\par\end{center}}
\preauthor{\begin{center}\normalsize\begin{tabular}{c}}
\postauthor{\end{tabular}\par\end{center}}

\predate{\begin{center}\normalsize}
\postdate{\par\end{center}}

\title{Saturation and Semi-Saturation Numbers
of Joins with Complete Graphs}

\author{
Changqing Xu \qquad Yan Yang \qquad Junxue Zhang* \\[0.5em]
{\small School of Science, Hebei University of  Technology, Tianjin 300401, China} 
}

\date{\today}

\begin{document}

\maketitle
\begin{sloppypar}
\baselineskip 0.65cm
\renewcommand{\thefootnote}{}

\footnotetext{
Corresponding author. 
\\
Email: 
chqxu@hebut.edu.cn (C. Xu), yy21101004@163.com (Y. Yang),  jxzhang@hebut.edu.cn (J. Zhang)}


\begin{abstract}
For a nonempty family $\mathcal F$ of graphs, a graph is $\mathcal F$-semi-saturated if adding any missing edge
creates a copy of some member of $\mathcal F$ containing the added
edge, and is $\mathcal F$-saturated if it is also $\mathcal F$-free.
The minimum numbers of edges in such graphs of order $n$ are denoted
by $\operatorname{ssat}(n,\mathcal F)$ and
$\operatorname{sat}(n,\mathcal F)$, respectively.
We study these parameters for families of graphs obtained by taking the join of a fixed complete graph with each member of a given graph family.
Assume that every member of $\mathcal F$ has at least $8$
vertices, minimum degree at least two, and no triangle containing
two vertices of degree two.
We give explicit sufficient conditions under which
$
\operatorname{sat}(n,K_t\vee\mathcal F)
=\binom{t}{2}+t(n-t)+\operatorname{sat}(n-t,\mathcal F),
$
and characterize all extremal graphs.
For the family $\mathcal C_{\ge r}=\{C_\ell:\ell\ge r\}$ of
subdivisions of $C_r$, this identity holds for fixed $t\ge1$ and
$r\ge4t^2+8t+29$ and all sufficiently large $n$.
For $K_t\vee C_k$ with $k\ge2t^2+3t+4$, we establish the
corresponding identity for saturation when $n\ge2(k+t-2)^3$
and for semi-saturation when $n\ge56k^3+t-1$.
\end{abstract}



\section{Introduction}
All graphs considered in this paper are finite, simple, and undirected.
Let $\mathcal F$ be a nonempty family of graphs.
A graph  is \emph{$\mathcal F$-free} if it contains no subgraph
isomorphic to any member of $\mathcal F$.
A graph $G$  is \emph{$\mathcal F$-semi-saturated} if adding any missing edge
$e$ to $G$ creates a copy of some member of $\mathcal F$ containing
$e$, and \emph{$\mathcal F$-saturated} if it is also
$\mathcal F$-free.
The \emph{semi-saturation number} $\operatorname{ssat}(n,\mathcal F)$
and the \emph{saturation number}
$\operatorname{sat}(n,\mathcal F)$ are the minimum numbers of
edges in $n$-vertex $\mathcal F$-semi-saturated and
$\mathcal F$-saturated graphs, respectively.
Since every $\mathcal F$-saturated graph is
$\mathcal F$-semi-saturated, we have
$
\operatorname{ssat}(n,\mathcal F)
\le\operatorname{sat}(n,\mathcal F).
$
We write $\operatorname{SAT}(n,\mathcal F)$ for the family of
all $\mathcal F$-saturated graphs of order $n$ with
$\operatorname{sat}(n,\mathcal F)$ edges.
When $\mathcal F=\{F\}$, we write $F$ in place of $\mathcal F$.

The \emph{join} $G\vee H$ of two vertex-disjoint graphs $G$ and $H$
is obtained from their disjoint union by adding all edges between
$V(G)$ and $V(H)$.
For a graph family $\mathcal F$, we write
$
K_t\vee\mathcal F=\{K_t\vee F:F\in\mathcal F\}.
$ A classical result of Erd\H{o}s, Hajnal, and Moon~\cite{EHM1964}
states that, for $n\ge r\ge3$,
$
\operatorname{sat}(n,K_r)
=(r-2)n-\binom{r-1}{2} 
$
with $K_{r-2}\vee\overline{K}_{n-r+2}$ as the unique
extremal graph, where $\overline{K}_{n-r+2}$ denotes the
empty graph on $n-r+2$ vertices.  Since $K_r=K_t\vee K_{r-t}$ for $1\le t<r$,
a natural extension is to study
$\operatorname{sat}(n,K_t\vee F)$ for a general graph $F$. 

Several special cases have been studied.
Chen, Faudree, and Gould~\cite{CFG} determined
the saturation numbers of books $K_2\vee\overline{K}_p$
and generalized books $K_t\vee\overline{K}_p$
for sufficiently large $n$.
Faudree, Ferrara, Gould, and Jacobson~\cite{FFGJ}
determined the saturation numbers of generalized
friendship graphs $K_t\vee rK_q$
for sufficiently large $n$, where $rK_q$ denotes
the disjoint union of $r$ copies of $K_q$.
Further results on joins of complete graphs
and disjoint unions of cliques were obtained
by Hu, Ji, and Zhang~\cite{HJZ}.

For a general graph $F$ with at least one edge,
repeated application of a result of Cameron and
Puleo~\cite{CP} gives
$
\operatorname{sat}(n,K_t\vee F)
\le \binom{t}{2}+t(n-t)+\operatorname{sat}(n-t,F)
$
for $t\ge1$ and $n\ge t+|V(F)|$.
It is natural to ask when equality holds.
Hu, Luo, and Peng~\cite{HLP} proved that if $F$ has no
isolated vertices and
$n\ge3t^2-t+2\operatorname{sat}(n-t,F)+1$, then
$
\operatorname{sat}(n,K_t\vee F)
=\binom{t}{2}+t(n-t)+\operatorname{sat}(n-t,F).
$
This condition applies when $\operatorname{sat}(n,F)$
is sufficiently small relative to $n$. 
For cycles, however, their condition fails when $n$
is sufficiently large.
Indeed, the following results hold for sufficiently
large $n$. 
 Ollmann~\cite{Ollmann} and Tuza~\cite{Tuza}
determined that
$
\operatorname{sat}(n,C_4)
=\left\lfloor\frac{3n-5}{2}\right\rfloor
$  
and Chen~\cite{Chen} proved that
$
\operatorname{sat}(n,C_5)
=\left\lceil\frac{10(n-1)}7\right\rceil
$.
For $C_6$, Lan, Shi, Wang, and Zhang~\cite{LSWZ}
showed that
$
\frac{4n}{3}-2
\le\operatorname{sat}(n,C_6)
\le\frac{4n+1}{3}
$.
For longer cycles, F\"uredi and Kim~\cite{FK} showed that, for
$k\ge7$ and $n\ge2k-5$,
\[
\left(1+\frac1{k+2}\right)n-1
<
\operatorname{sat}(n,C_k)
<
\left(1+\frac1{k-4}\right)n+\binom{k-4}{2}.
\]
The lower bound makes the condition above impossible
for $F=C_k$ when $n$ is sufficiently large. 
More recently, Mohammadian, Poursoltani, and
Tayfeh-Rezaie~\cite{MPT} proved that, for each fixed
even integer $k\ge28$,
$
\operatorname{sat}(n,C_k)
=\frac{k-3}{k-4}n+O(1).
$

When $t=1$, the join $K_t\vee C_k$ is the wheel
$W_k=K_1\vee C_k$.
Song, Hu, Ji, and Cui~\cite{SHJC} determined
$\operatorname{sat}(n,W_4)$ and characterized
all extremal graphs.
For longer wheels, Qiu, He, Lu, and Xu~\cite{QHLX}
proved that, for $k\ge8$ and $n\ge56k^3$,
$
\operatorname{sat}(n,W_k)
=n-1+\operatorname{sat}(n-1,C_k)
$
and
$
\operatorname{ssat}(n,W_k)
=n-1+\operatorname{ssat}(n-1,C_k).
$
They also showed that every extremal graph for either
parameter contains a universal vertex.
A vertex of a graph is called \emph{universal} if it is
adjacent to every other vertex of the graph.
These results lead to the question of whether
the wheel identities extend to $K_t\vee C_k$
for $t\ge2$.

Our main result establishes a general saturation identity and characterizes all extremal graphs for joins of complete graphs with graph families whose members satisfy suitable minimum degree and triangle conditions (Theorem~\ref{mainthm:Kt+H}).
We apply this result to subdivisions of cycles and to fixed cycles in Theorems~\ref{mainthm:Kt+Cger} and~\ref{mainthm:Kt+Ck}, respectively.


\begin{thm}\label{mainthm:Kt+H}
Let $t,f,n$ be positive integers with $f\ge8$, and let
$\mathcal F$ be a nonempty family of graphs.
Assume that every $F\in\mathcal F$ has at least $f$ vertices,
satisfies $\delta(F)\ge2$, and contains no triangle with two
vertices of degree two in $F$.
If $n\geq2(f+t-2)^3$ and
\(
\operatorname{sat}(n-t,\mathcal F)
<
\left(\frac32-\frac{t^2+2t-1}{f+t-2}\right)n,
\)
then
\[
\operatorname{sat}(n,K_t\vee \mathcal F)
=
\binom{t}{2}+t(n-t)+\operatorname{sat}(n-t,\mathcal F).
\]
Moreover,
\(\operatorname{SAT}(n,K_t\vee \mathcal F)
=
\bigl\{K_t\vee G_1:
G_1\in\operatorname{SAT}(n-t,\mathcal F)\bigr\}.\)

\end{thm}

We next consider joins of complete graphs with subdivisions of cycles.
For an integer $r\ge3$, let
$\mathcal C_{\ge r}=\{C_\ell:\ell\ge r\}$,
which is precisely the family of all subdivisions of $C_r$.
Ferrara et al.~\cite{FJMTW} proved that, for every fixed
integer $r\ge36$ and all sufficiently large $n$,
$ 
\operatorname{sat}(n,\mathcal C_{\ge r})
\le
\left(\frac54+\frac8r\right)n.
 $
Combining this bound with Theorem~\ref{mainthm:Kt+H} gives the following result.

\begin{thm}\label{mainthm:Kt+Cger}
Let $t$ and $r$ be fixed positive integers with $r\ge4t^2+8t+29$, and let
$\mathcal C_{\ge r}=\{C_\ell:\ell\ge r\}$.
Then, for all sufficiently large $n$,
 $
\operatorname{sat}(n,K_t\vee\mathcal C_{\ge r})
=
\binom{t}{2}+t(n-t)
+\operatorname{sat}(n-t,\mathcal C_{\ge r}).
 $
Moreover,
 $
\operatorname{SAT}(n,K_t\vee\mathcal C_{\ge r})
=
\bigl\{
K_t\vee G_1:
G_1\in\operatorname{SAT}(n-t,\mathcal C_{\ge r})
\bigr\}.
 $
\end{thm}

For joins with a fixed cycle, we obtain explicit conditions on
$t$, $k$, and $n$ for both saturation and semi-saturation.

\begin{thm}\label{mainthm:Kt+Ck}
Let $t,k,n$ be positive integers with $k\ge2t^2+3t+4$.
Then the following hold. \begin{enumerate}
\item[\rm{(i)}] If $n\ge2(k+t-2)^3$, then
$
\operatorname{sat}(n,K_t\vee C_k)
=\binom{t}{2}+t(n-t)+\operatorname{sat}(n-t,C_k).
$
Moreover,
$
\operatorname{SAT}(n,K_t\vee C_k)
=\{K_t\vee G_1:G_1\in\operatorname{SAT}(n-t,C_k)\}.
$

\item[\rm{(ii)}] If $n\ge56k^3+t-1$, then
$
\operatorname{ssat}(n,K_t\vee C_k)
=\binom{t}{2}+t(n-t)+\operatorname{ssat}(n-t,C_k).
$
\end{enumerate}
\end{thm}

Combining Theorem~\ref{mainthm:Kt+Ck}\textup{(i)} with the result of Mohammadian,
Poursoltani, and Tayfeh-Rezaie~\cite{MPT} gives the
following corollary. 

\begin{cor}
Let $t$ be a fixed integer with $t\ge1$, and let
$k$ be a fixed even integer with $k\ge\max\{28,2t^2+3t+4\}$.
Then, for all sufficiently large $n$, 
\[
\operatorname{sat}(n,K_t\vee C_k)
=
\left(t+\frac{k-3}{k-4}\right)n+O(1),
\]
where the constant in the $O(1)$ term depends only on $t$ and $k$.
\end{cor}

The paper is organized as follows.
Section~2 states the general lower bound in
semi-saturated graphs with fewer than $t$ universal vertices and gives reduction formulas for saturated and semi-saturated graphs
with at least $t$ universal vertices.
Section~3 applies these results to prove
Theorems~\ref{mainthm:Kt+H}--\ref{mainthm:Kt+Ck}.
Section~4 establishes the structural lemmas and proves the
general lower bound stated in Section~2.

\section{Preliminaries}
We first state a lower bound for semi-saturated graphs with fewer than $t$ universal vertices. Its proof is deferred to Section~4.
\begin{lem}\label{thm:K_t+Hlower}
Let $t,f,n$ be positive integers with $f\ge8$, and let
$\mathcal F$ be a nonempty family of graphs.
Assume that every $F\in\mathcal F$ has at least $f$ vertices,
satisfies $\delta(F)\ge2$, and contains no triangle with two
vertices of degree two in $F$.
If $n\geq2(f+t-2)^3$ and $G$ is an $n$-vertex $(K_t\vee \mathcal F)$-semi-saturated graph
with fewer than $t$ universal vertices, then
\[
e(G)\ge
\left(t+\frac32-\frac{t^2+2t-1}{f+t-2}\right)n.
\]
\end{lem} 

We next consider graphs with at least $t$ universal vertices.
For a single forbidden graph $F$, Cameron and Puleo~\cite{CP}
established the upper bound
$
\operatorname{sat}(n,K_1\vee F)
\le n-1+\operatorname{sat}(n-1,F).
$
Hu, Ji, and Zhang~\cite{HJZ} proved that, if $v$ is a universal
vertex of $G$, then $G$ is $(K_1\vee F)$-saturated if and only
if $G-v$ is $F$-saturated.
 We give a proof for graph families and  general $t$.
For $X \subseteq V(G)$, let $G[X]$ be the subgraph of $G$ whose vertex set is $X$ and whose edge set consists of all edges of $G$ which have both ends in $X$, and let $G-X=G[V(G)\setminus X]$. For any graph $F$, in each copy of $K_t\vee F$
considered below, we fix a $K_t$ part and the corresponding
$F$ part, and call
the vertices in the former its \emph{centers}.
Each center has degree $t+|V(F)|-1$ in that copy.

\begin{lem}
\label{lem:universal-join}
Let $t,n$ be positive integers, and let $\mathcal H$ be a
nonempty family of graphs, each with at least one edge.
If $n\geq t+\min\{|V(H)|:H\in\mathcal H\}$, then the following hold.
\begin{enumerate}
\item[\textup{(i)}]
$\operatorname{sat}(n,K_t\vee\mathcal H)
\leq\binom{t}{2}+t(n-t)+\operatorname{sat}(n-t,\mathcal H)$.

\item[\textup{(ii)}]
The minimum number of edges in a
$(K_t\vee\mathcal H)$-saturated graph of order $n$
with at least $t$ universal vertices is
$\binom{t}{2}+t(n-t)+\operatorname{sat}(n-t,\mathcal H)$.
Moreover, such a graph has the minimum number of edges
if and only if it is isomorphic to $K_t\vee G_1$ for some
$G_1\in\operatorname{SAT}(n-t,\mathcal H)$.
\end{enumerate}
\end{lem}

\begin{proof}
Let $G$ be a graph of order $n$ with at least $t$ universal
vertices. Choose a set $U$ of $t$ universal vertices of $G$,
and let $G_0$ denote $G-U$.
We first show that $G$ is $(K_t\vee\mathcal H)$-free
if and only if $G_0$ is $\mathcal H$-free.
If $G_0$ contains a copy of some $H\in\mathcal H$,
the join of this copy and $G[U]$ gives a copy of $K_t\vee H$ in $G$.
Conversely, suppose that $G$ contains a copy $T$ of
$K_t\vee H$ for some $H\in\mathcal H$.
Let $T_1,T_2$ denote its $K_t$ part and $H$ part,
respectively, and let $A$ denote $U\cap V(T_2)$.
Since  $T_1,T_2$ are vertex-disjoint and
$|V(T_1)|=|U|=t$, we have
$|A|\leq t-|U\cap V(T_1)|
=|V(T_1)\setminus U|$.
Replace the vertices of $A$ in $T_2$ by distinct
vertices of $V(T_1)\setminus U$.
Since $T_1$ is complete and every vertex of $T_1$
is adjacent to every vertex of $T_2$, this gives
a copy of $H$ in $G_0$. This proves the stated equivalence.

Now assume that $G$ is $(K_t\vee\mathcal H)$-free,
so that $G_0$ is $\mathcal H$-free.
As every vertex of $U$ is universal,
$E(\overline G)=E(\overline{G_0})$.
For each $e\in E(\overline{G_0})$, the vertices of $U$
remain universal in $G+e$.
Applying the above argument to $G+e$, we see that $G+e$ contains a copy of
a member of $K_t\vee\mathcal H$ if and only if $G_0+e$
contains a copy of a member of $\mathcal H$.
Hence, $G$ is $(K_t\vee\mathcal H)$-saturated
if and only if $G_0$ is $\mathcal H$-saturated.

For \textup{(i)}, choose
$G_1\in\operatorname{SAT}(n-t,\mathcal H)$.
By the saturation equivalence proved above, $K_t\vee G_1$ is
$(K_t\vee\mathcal H)$-saturated and has at least
$t$ universal vertices.
Its number of edges is
$\binom{t}{2}+t(n-t)+\operatorname{sat}(n-t,\mathcal H)$.
This proves \textup{(i)} and gives the upper bound
in \textup{(ii)}.

For the lower bound in \textup{(ii)}, if $G$ is
$(K_t\vee\mathcal H)$-saturated, then $G_0$ is
$\mathcal H$-saturated.
Therefore,
$e(G)=\binom{t}{2}+t(n-t)+e(G_0)
\geq\binom{t}{2}+t(n-t)+\operatorname{sat}(n-t,\mathcal H)$,
with equality if and only if
$G_0\in\operatorname{SAT}(n-t,\mathcal H)$.
\end{proof}

For cycles, we obtain the corresponding reduction for semi-saturation using the wheel theorem.

\begin{lem}
\label{lemma:universal}
Let $t,k,n$ be positive integers.
If $k\geq8$ and $n\geq56k^3+t-1$, then the following hold.
\begin{enumerate}
\item[\textup{(i)}]
$\displaystyle
\operatorname{ssat}(n,K_t\vee C_k)
\le \binom{t}{2}+t(n-t)+\operatorname{ssat}(n-t,C_k).
$

\item[\textup{(ii)}]
The minimum number of edges in a
$(K_t\vee C_k)$-semi-saturated graph of order $n$
with at least $t$ universal vertices is
$
\binom{t}{2}+t(n-t)+\operatorname{ssat}(n-t,C_k).
$
\end{enumerate}
\end{lem}
\begin{proof}
For (i),   choose a $C_k$-semi-saturated graph
$G_0$ of order $n-t$ with
$e(G_0)=\operatorname{ssat}(n-t,C_k)$.
For every $e\in E(\overline{K_t\vee G_0})$, we have
$e\in E(\overline{G_0})$, so $G_0+e$ contains
a copy of $C_k$ using $e$.
Together with the added $K_t$, this gives a copy
of $K_t\vee C_k$ using $e$ in $(K_t\vee G_0)+e$.
Thus, $\operatorname{ssat}(n,K_t\vee C_k)\leq\binom{t}{2}+t(n-t)+\operatorname{ssat}(n-t,C_k)$. This proves \textup{(i)} and gives the upper bound
in \textup{(ii)}.

For the lower bound in \textup{(ii)}, we first claim
that, for $j\geq2$, if $G$ is
$(K_j\vee C_k)$-semi-saturated and $u$ is a universal
vertex of $G$, then $G-u$ is
$(K_{j-1}\vee C_k)$-semi-saturated.
Consider $xy\in E(\overline{G-u})$.
Since $G$ is $(K_j\vee C_k)$-semi-saturated, $G+xy$ contains a copy $T$
of $K_j\vee C_k$ using $xy$.
If $u$ is not a center, choose a center $h$ of $T$ as follows. Note that  $j\geq2$.
If exactly one of $x,y$ is a center, choose $h$ different
from that endpoint; if both are centers, choose
$h\in\{x,y\}$; otherwise, choose any center.
If $u$ lies on the cycle of $T$, replace $u$ on
the cycle by $h$, and replace the center $h$ by $u$.
If $u\notin V(T)$, replace a cycle vertex different
from $x,y$ by $h$, and replace the center $h$ by $u$.
Since $u$ is universal and $h$ is adjacent to every
vertex of the cycle of $T$, this gives another
copy of $K_j\vee C_k$.
If $xy$ is a cycle edge, neither endpoint is replaced;
otherwise, both endpoints are retained and at least
one remains a center.
Thus, the new copy still uses $xy$ and has $u$
as a center.
Deleting $u$ gives a copy of $K_{j-1}\vee C_k$
using $xy$ in $(G-u)+xy$, proving the claim.

Now let $G$ be a $(K_t\vee C_k)$-semi-saturated graph
of order $n$ with at least $t$ universal vertices.
Delete $t-1$ of these vertices successively,
and denote the resulting graph by $G_1$,
taking $G_1=G$ when $t=1$.
The remaining selected vertices are still universal
after each deletion.
Repeated application of the claim shows that $G_1$
is a $W_k$-semi-saturated graph of order $n-t+1$,
where $W_k=K_1\vee C_k$.
Since $k\geq8$ and $n-t+1\geq56k^3$,
the wheel theorem of Qiu, He, Lu, and Xu~\cite{QHLX}
gives $e(G_1)\geq\operatorname{ssat}(n-t+1,W_k)
=n-t+\operatorname{ssat}(n-t,C_k)$.
Consequently,
$e(G)=e(G_1)+\binom{t-1}{2}+(t-1)(n-t+1)
\geq\binom{t}{2}+t(n-t)+\operatorname{ssat}(n-t,C_k)$,
completing the proof.
\end{proof}

\section{Proofs of the main results}
 We combine Lemmas~\ref{thm:K_t+Hlower} and~\ref{lem:universal-join} to prove Theorem~\ref{mainthm:Kt+H}.
\subsection{Proof of Theorem~\ref{mainthm:Kt+H}}

\begin{proof}

By Lemma~\ref{lem:universal-join} \textup{(i)}, we have
$\operatorname{sat}(n,K_t\vee\mathcal F)
\leq\binom{t}{2}+t(n-t)+\operatorname{sat}(n-t,\mathcal F)$.
Using
$\operatorname{sat}(n-t,\mathcal F)
<(\frac32-\frac{t^2+2t-1}{t+f-2})n$,
we obtain
\[
\begin{aligned}
\operatorname{sat}(n,K_t\vee \mathcal F)
&<
\left(t+\frac32-\frac{t^2+2t-1}{t+f-2}\right)n
-\frac{t(t+1)}2\\
&<
\left(t+\frac32-\frac{t^2+2t-1}{t+f-2}\right)n.
\end{aligned}
\]

Let $G\in\operatorname{SAT}(n,K_t\vee \mathcal F)$.
Suppose that $G$ has fewer than
$t$ universal vertices.
The graph $G$ is also
$(K_t\vee\mathcal F)$-semi-saturated.
By Lemma~\ref{thm:K_t+Hlower},  we have
$e(G)\geq
(t+\frac32-\frac{t^2+2t-1}{t+f-2})n$,
contradicting the strict upper bound above.
Thus, $G$ has at least $t$ universal vertices.
By Lemma~\ref{lem:universal-join} \textup{(ii)},
\[
e(G)=\operatorname{sat}(n,K_t\vee \mathcal F)
=
\binom{t}{2}+t(n-t)+\operatorname{sat}(n-t,\mathcal F),
\]
and  $G \cong K_t\vee G_1$  for some
$G_1\in\operatorname{SAT}(n-t,\mathcal F)$.
\end{proof}

\subsection{Proof of Theorem~\ref{mainthm:Kt+Cger}}
\begin{proof}
    
Every member of $\mathcal C_{\geq r}$ is a triangle-free graph
with at least $r$ vertices and minimum degree two.
Since $r\geq4t^2+8t+29$, we have
\[
\frac{t^2+2t-1}{r+t-2}+\frac8r
<
 \frac14.
\]

For all sufficiently large $n$, we have
$n\geq2(r+t-2)^3$, and the upper bound of
Ferrara et al.~\cite{FJMTW} gives
\[
\operatorname{sat}(n-t,\mathcal C_{\geq r})
\leq
\left(\frac54+\frac8r\right)(n-t)
<
\left(\frac32-\frac{t^2+2t-1}{r+t-2}\right)n.
\]
By Theorem~\ref{mainthm:Kt+H}, we obtain
the stated saturation identity and extremal graphs.
\end{proof}

\subsection{Proof of Theorem~\ref{mainthm:Kt+Ck}}
We apply Theorem~\ref{mainthm:Kt+H} to prove part~(i), and use
Lemmas \ref{thm:K_t+Hlower} and   \ref{lemma:universal} to prove part~(ii). 

\begin{proof}
We first prove (i).  Since $k\geq2t^2+3t+4$ and $n\geq2(k+t-2)^3$,
we have
\[
\left(\frac12-\frac{t^2+2t-1}{k+t-2}
-\frac1{k-4}\right)n
>
\frac{n}{4(k+t-2)}
\geq\frac{(k+t-2)^2}{2}
>\binom{k-4}{2}.
\]
By the upper bound of F\"uredi and Kim~\cite{FK},
we obtain
\[
\operatorname{sat}(n-t,C_k)
<
\left(1+\frac1{k-4}\right)(n-t)+\binom{k-4}{2}
<
\left(\frac32-\frac{t^2+2t-1}{k+t-2}\right)n.
\]
The cycle $C_k$ has minimum degree two and contains
no triangle. By Theorem~\ref{mainthm:Kt+H}
with   the singleton family $\mathcal F=\{C_k\}$, we obtain $\operatorname{sat}(n,K_t\vee C_k)
=
\binom{t}{2}+t(n-t)+\operatorname{sat}(n-t,C_k)$ and $\operatorname{SAT}(n,K_t\vee C_k)
=
\bigl\{K_t\vee G_1:
G_1\in\operatorname{SAT}(n-t,C_k)\bigr\}.$

We next prove (ii).
Note that $k\geq8$ and $n\geq56k^3+t-1\ge2(k+t-2)^3$.
By Lemma~\ref{lemma:universal} \textup{(i)}, we have
\[
\begin{aligned}
\operatorname{ssat}(n,K_t\vee C_k)
&\leq\binom{t}{2}+t(n-t)
+\operatorname{ssat}(n-t,C_k)\\
&\leq\binom{t}{2}+t(n-t)
+\operatorname{sat}(n-t,C_k)\\
&<
\left(t+\frac32-\frac{t^2+2t-1}{k+t-2}\right)n.
\end{aligned}
\]
Let $G$ be a $(K_t\vee C_k)$-semi-saturated graph
of order $n$ such that
$e(G)=\operatorname{ssat}(n,K_t\vee C_k)$.
Assume that $G$ has fewer than
$t$ universal vertices.
By Lemma~\ref{thm:K_t+Hlower}, we have
$e(G)\geq
(t+\frac32-\frac{t^2+2t-1}{k+t-2})n$,
contradicting the strict upper bound above.
Thus, $G$ has at least $t$ universal vertices.
By Lemma~\ref{lemma:universal} \textup{(ii)}, we obtain
$
e(G)=\operatorname{ssat}(n,K_t\vee C_k)
=\binom{t}{2}+t(n-t)+\operatorname{ssat}(n-t,C_k).
$
\end{proof}

\section{Proof of the general lower bound}
We prove Lemma~\ref{thm:K_t+Hlower} in this section.
 We first establish structural
properties of low-degree vertices and their neighborhoods,
and then use these properties to derive the required
edge bound.

\subsection{Structural lemmas}
For $v\in V(G)$, define
$N_G(v)=\{u\in V(G):uv\in E(G)\}$. 
Let $X$ and $Y$ be sets of vertices of a graph $G$. 
We define
\(
N_G(X)=\left(\bigcup_{v\in X}N_G(v)\right)\setminus X\) and \(
 N_G[X]=X\cup N_G(X).
\)
We denote by $E_G[X,Y]$ the set of edges of $G$ with one end in $X$ and the other end in $Y$, and by $e_G(X,Y)$ their number.
Let $x$ and $y$ be two vertices of $G$. The \emph{distance}
between $x$ and $y$, denoted by $d_G(x,y)$, is the length
of a shortest $xy$-path; if no such path exists, we set
$d_G(x,y)=\infty$.
Throughout the remainder of this section, let $t$, $f$, and $n$ be positive integers such that $f\geq8$ and $n\geq 2(f+t-2)^3$. Let $\mathcal F$ be a nonempty family of graphs such that
every $F\in\mathcal F$ has at least $f$ vertices,
satisfies $\delta(F)\geq2$, and contains no triangle
with two vertices of degree two in $F$.
Let $G$ be a $(K_t\vee\mathcal F)$-semi-saturated graph
of order $n$ with fewer than $t$ universal vertices.

We first observe that $\delta(G)\geq t+1$.
If a vertex $v$ is not universal, choose a vertex $u$
nonadjacent to $v$.
Since $G$ is $(K_t\vee\mathcal F)$-semi-saturated,
$G+uv$ contains a copy $T$ of $K_t\vee F$ using $uv$,
for some $F\in\mathcal F$.
As $\delta(F)\geq2$, we have $d_G(v)+1\geq d_T(v)\geq t+2$.
If $v$ is universal, then
$d_G(v)=n-1\geq t+f-1$.
We partition $V(G)$ into the following three sets:
$L_G=\{v\in V(G):d_G(v)=t+1\},
M_G=\{v\in V(G):t+2\leq d_G(v)\leq t+f-3\},
H_G=\{v\in V(G):d_G(v)\geq t+f-2\}.$
Let
\(
H_G^*=L_G\cup M_G,
\)
and let
\(
D=\{v\in V(G):d_G(v)=n-1\}
\).

\begin{lem}{\rm{\label{lem1}}}
\textit{For two distinct vertices \(u\) and \(v\) of \(G\), we have \(d_G(u,v)\leq 2\).}
\end{lem}
\begin{proof}
Let $u$ and $v$ be distinct vertices of $G$.
If they are adjacent, then $d_G(u,v)=1$.
Otherwise, since $G$ is $(K_t\vee \mathcal F)$-semi-saturated, $G+uv$ contains a copy
of $K_t\vee F$ using $uv$ for some $F\in\mathcal F$.
Since $\delta(F)\geq2$, $u$ and $v$
have at least $t$ common neighbors in $G$.
As $t\geq1$, it follows that $d_G(u,v)=2$.
\end{proof}

\begin{lem}{\rm{\label{lem2}}}
The set $L_G$ is independent, and any two distinct
vertices in $L_G$ have exactly $t$ common neighbors.
\end{lem}
\begin{proof}
Suppose, to the contrary, that two distinct vertices $u$ and $v$ in $L_G$ are adjacent. Since $d_G(u)=t+1$ and
$n\geq t+f$, there exists a vertex $w$ outside
$N_G[u]$. Since $G$ is $(K_t\vee\mathcal F)$-semi-saturated, $G+uw$ contains
a copy $T_1$ of $K_t\vee F_1$ using $uw$ for some $F_1\in\mathcal F$.
Since $\delta(T_1)\geq t+2$ and $d_{G+uw}(u)= t+2$, all edges
incident with $u$ in $G+uw$ belong to $T_1$.
In particular, $v\in V(T_1)$. As $v$ is not
incident with $uw$, we have
$d_G(v)\geq d_{T_1}(v)\geq t+2$, a contradiction.
Thus, $L_G$ is independent.

Let $u$ and $v$ be distinct vertices in $L_G$. 
Since $G$ is $(K_t\vee\mathcal F)$-semi-saturated, $G+uv$ contains a copy $T_2$
of $K_t\vee F_2$ using $uv$ for some $F_2\in\mathcal F$.
Since both vertices have degree $t+2$ in $G+uv$
and $f\ge 8$, they belong to the $F_2$ part of $T_2$.
The $t$ centers of $T_2$ are common neighbors
of $u$ and $v$ in $G$, so
$|N_G(u)\cap N_G(v)|\geq t$. Note that both vertices have
degree $t+1$ in $G$. If this intersection has more than $t$ vertices,
then $N_G(u)=N_G(v)$. Their common neighborhood
consists of the centers and one additional
vertex $x$. Thus, within the $F_2$ part of $T_2$, the
neighbors of $u$ are $v$ and $x$, and those
of $v$ are $u$ and $x$. Hence the $F_2$ part contains
a triangle $xuv$ in which both $u$ and $v$ have degree
two, contrary to the assumption on $F_2$.
Therefore, $|N_G(u)\cap N_G(v)|=t$.
\end{proof}

\begin{lem}{\rm{\label{lem4}}}
Assume that $|D|<t$ and $H_G^*\neq\varnothing$.
If
\(
\left(\bigcap_{v\in H_G^*}N_G(v)\right)
\setminus D=\varnothing,
\)
then there exists a set $X\subseteq V(G)$ such that
$|X|\leq(t+f-2)^2$ and
$|N_G(v)\cap X|\geq t+1$ for every   $v\in V(G) \setminus X$.
\end{lem}
\begin{proof}
Choose \(x\in H_G^*\), and write
\((N_G(x)\cap H_G)\setminus D=\{x_1,x_2,\ldots,x_r\}\),
where \(r\le f+t-3\). Since
\(\Bigl(\bigcap_{v\in H_G^*}N_G(v)\Bigr)\setminus D=\varnothing,\)
for each \(i\in \{1,2,\dots, r\}\), there exists \(y_i\in H_G^*\)
such that \(x_i y_i\notin E(G)\). Set
\(X=N_G[x]\cup(\bigcup_{i=1}^{r}N_G[y_i]).\)
By the degree bounds for the vertices in \(H_G^*\), we have
\(|X|\le (f+t-2)+(f+t-3)(f+t-2)=(f+t-2)^2.\)

Next, we prove that $|N_G(v)\cap X|\geq t+1$ for all \(v\in V(G)\setminus X\).
Let \(v\in V(G)\setminus X\). Since \(v\notin X\), we have
\(v\notin N_G[x]\). Since \(G\) is \((K_t\vee\mathcal F)\)-semi-saturated,
\(G+vx\) contains a copy \(T_1\) of \(K_t\vee F_1\) using \(vx\) for some $F_1\in\mathcal F$. The degree of $x$ in $G+vx$ is at most $f+t-2$, so the vertex $x$ belongs to the $F_1$ part of $T_1$. If \(v\) is a center of \(T_1\), then
\(|N_G(x)\cap N_G(v)|\ge t+1\), as required. Otherwise, \(|N_G(x)\cap N_G(v)|\ge t\). Since
\(\{x_1,x_2,\ldots,x_r\}\cup D\subseteq N_G(x)\), where \(|D|<t\), there exists
some vertex among \(x_1,x_2,\ldots,x_r\) that is a center of \(T_1\);
without loss of generality, let \(x_1\) be such a vertex. Since
\(v\notin N_G[y_1]\), \(G+vy_1\) contains a copy \(T_2\) of
\(K_t\vee F_2\) using \(vy_1\) for some $F_2\in\mathcal F$. The degree of $y_1$ in $G+vy_1$ is at most $f+t-2$, so $y_1$ belongs to the $F_2$ part of $T_2$. If $v$ is a center of $T_2$, then \(|N_G(y_1)\cap N_G(v)|\ge t+1\), as desired. Otherwise, all centers of $T_2$ belong to
$N_G(y_1)\cap N_G(v)$.
Since \(x_1y_1\notin E(G)\), it follows that \(x_1\) is not a center of
\(T_2\). The union of the center sets of \(T_1\) and \(T_2\) contains at least \(t+1\) vertices. Therefore, $|N_G(v)\cap X|\geq t+1$.
\end{proof}

\begin{lem}{\rm{\label{lem5}}}
Assume that $|D|<t$ and $|H_G^*|>f+t-1$.
For $x\in L_G$ with $N_G(x)=\{x_1,x_2,\ldots,x_{t+1}\}$, we have
$d_G(x_i)\geq f+t-1$ for $1\leq i\leq t+1$, and $x_ix_j\in E(G)$ for
all distinct $i,j\in\{1,2,\ldots,t+1\}$.
\end{lem}
\begin{proof}
Let \(x\in L_G\) and \(N_G(x)=\{x_1,x_2,\ldots,x_{t+1}\}\). Since
\(|H_G^*|>f+t-1\), choose $y\in H_G^*\setminus N_G[x]$. Hence \(xy\notin E(G)\). Since \(G\) is
\((K_t\vee\mathcal F)\)-semi-saturated, \(G+xy\) contains a copy $T_1$ of
\(K_t\vee F_1\) using \(xy\) for some $F_1\in \mathcal F$. 
Since $d_{G+xy}(x)=t+2<t+f-1$, the vertex $x$
belongs to the $F_1$ part of $T_1$.
Moreover, $\delta(T_1)\geq t+2$, so all neighbors
of $x$ in $G+xy$ belong to $T_1$. Note that $d_{G+xy}(y)\leq f+t-2$. 
Thus, \(t\) vertices among \(x_1,x_2,\ldots,x_{t+1}\) are centers of
\(T_1\); without loss of generality, let them be \(x_1,x_2,\ldots,x_t\).
It follows that \(d_G(x_i)\ge f+t-1\) for \(1\le i\le t\), and
\(x_i x_j\in E(G)\) for all distinct \(i,j\in\{1,2,\ldots,t+1\}\).

Since \(|D|<t\), there is a vertex among \(x_1,x_2,\ldots,x_t\)
that does not belong to \(D\). Without loss of generality, let
\(x_1\notin D\). Then there exists a vertex
\(z\) such that \(x_1z\notin E(G)\). Since each vertex in \(N_G[x]\setminus \{x_1\}\) is
adjacent to \(x_1\), we have \(z\notin N_G[x]\). Since \(G\) is
\((K_t\vee \mathcal F)\)-semi-saturated, \(G+xz\) contains a copy $T_2$ of
\(K_t\vee F_2\) using \(xz\) for some $F_2\in \mathcal F$. Since $d_{G+xz}(x)=t+2$, the vertex $x$ belongs to
the $F_2$ part of $T_2$  and all its neighbors belong to $T_2$. 
As $x_1z\notin E(G)$, neither $x_1$ nor $z$ can be
a center of $T_2$.
Thus, the $t$ centers of $T_2$ are precisely
$x_2,\ldots,x_{t+1}$.
Therefore,
\(d_G(x_{t+1})\ge f+t-1\). 
\end{proof}

\begin{lem}{\rm{\label{lem6}}}
There exists a set $L_G^*\subseteq N_G(L_G)$ with $|L_G^*|=|L_G|$. Moreover, if $w$ is adjacent to every vertex of $L_G$,
then $L_G^*$ can be chosen so that $w\notin L_G^*$.
\end{lem}
\begin{proof}
If $L_G=\varnothing$, take $L_G^*=\varnothing$. If $|L_G|=1$, choose a neighbor of its unique vertex
different from $w$. 
Now assume that $L_G$ consists of the vertices $x_1,x_2,\ldots,x_r$, where $r\geq2$, and write $A_i=N_G(x_i)$. By the definition of $L_G$ and Lemma~\ref{lem2}, each $A_i$ has size $t+1$, and $|A_i\cap A_j|=t$ whenever $i\neq j$. Let $Y$ denote $A_1\cap A_2$. Then $A_1=Y\cup\{z_1\}$ and $A_2=Y\cup\{z_2\}$, where $z_1$ and $z_2$ are distinct vertices outside $Y$. The intersection conditions imply that each $A_i$ either contains $Y$ or has the form $(Y\setminus\{y\})\cup\{z_1,z_2\}$ for some $y\in Y$. If every $A_i$ contains $Y$, write $A_i=Y\cup\{z_i\}$. The vertices $z_1,z_2,\ldots,z_r$ are pairwise distinct, so we may take $L_G^*=\{z_1,z_2,\ldots,z_r\}$. The common neighbor $w$ of the vertices in $L_G$
belongs to $Y$, whereas $L_G^*\cap Y=\varnothing$.
Hence $w\notin L_G^*$.
Otherwise, without loss of generality, assume that $A_3=(Y\setminus\{y_1\})\cup\{z_1,z_2\}$. 
Let $U=Y\cup\{z_1,z_2\}$. 
Each $A_i$ is a subset of $U$.
Otherwise, some $A_i=Y \cup \{z\}$ for some $z \notin Y\cup\{z_1, z_2\}$ and $|A_i\cap A_3|=|Y\setminus \{y_1\}|=t-1$,  a contradiction. 
  Since these sets are distinct and $|U|=t+2$, we have $3\leq r\leq t+2$. For each $i$, let $u_i$ denote the unique vertex in $U\setminus A_i$. Take $L_G^*=\{u_1,u_2,\ldots,u_r\}$. Since $u_i\notin A_i$ for each $i$, no vertex of $L_G^*$
is a common neighbor of the vertices in $L_G$.
Hence $w\notin L_G^*$.
The vertices $u_1,u_2,\ldots,u_r$ are pairwise distinct, and $L_G^*\subseteq U=A_1\cup A_2$. Therefore, $L_G^*\subseteq\bigcup_{i=1}^{r}A_i=N_G(L_G)$ and $|L_G^*|=|L_G|$.
\end{proof}

\begin{lem}{\rm{\label{lem3}}}
Assume that $|D|<t$ and $|H_G^*|>f+t-1$. For $R\subseteq V(G)$, we have
\(\sum_{v\in R}d_G(v)\ge (t+2)|R|+\max\left\{-|L_G|,\,-(n-|R|)\right\}.
\)
\end{lem}
\begin{proof}
Since each vertex $v\in V(G)\setminus L_G$ satisfies $d_G(v)\ge t+2$, we have $\sum_{v\in R}d_G(v)\ge (t+2)|R|-|R\cap L_G|$. 
By Lemma~\ref{lem5}, each vertex in $N_G(L_G)$ has degree at least $f+t-1$. By Lemma \ref{lem6}, there exists a set
$L_G^*\subseteq N_G(L_G)$ with $|L_G^*|=|L_G|$.
Since $|L_G^*\cap R|\ge |L_G^*|-|V(G)\setminus R|=|L_G|-(n-|R|)$, there are at least $\max\{0,|L_G|-(n-|R|)\}$ vertices in $R$ of degree at least $f+t-1$. 
Hence
\(\sum_{v\in R}d_G(v)\ge (t+2)|R|+(f-3)\max\left\{0,\ |L_G|-(n-|R|)\right\}-|L_G|.\)
Since $f>7$, \(\sum_{v\in R}d_G(v)\ge (t+2)|R|+\max\left\{-|L_G|,\,-(n-|R|)\right\}. 
\)
\end{proof}

\subsection{Proof of Lemma~\ref{thm:K_t+Hlower}}

\begin{proof}
Let $G$ be a $(K_t\vee\mathcal F)$-semi-saturated graph of
order $n$ with fewer than $t$ universal vertices.
We consider three cases according to the size of $H_G^*$
and whether a nonuniversal vertex is adjacent to every
vertex of $H_G^*$ to prove the following lower bound:
\begin{equation}
\label{eq:lower-bound}
e(G)\geq
\left(t+\frac32-\frac{t^2+2t-1}{t+f-2}\right)n.
\tag{1}
\end{equation}

\medskip
     \noindent\textbf{Case 1.} $|H_G^*|\leq f+t-1$.

Since $|H_G^*|\leq f+t-1$, we have
\(
2e(G)
\geq(t+1)|H_G^*|+(f+t-2)(n-|H_G^*|)
\geq(f+t-2)n-(f-3)(f+t-1).
\)
For $f\geq 8$ and $t\geq1$, we have
\[
\frac{f+t-2}{2}
-\left(t+\frac32-\frac{t^2+2t-1}{f+t-2}\right)
=
\frac{f^2-7f+t^2+t+8}{2(f+t-2)}
\geq\frac{f-3}{2(f+t-2)}.
\]
Consequently,
\[
e(G)\geq
\left(t+\frac32-\frac{t^2+2t-1}{f+t-2}\right)n
+\frac{f-3}{2}
\left(\frac{n}{f+t-2}-(f+t-1)\right).
\]
Since $n\geq2(f+t-2)^3$, the last term is positive.
Thus, inequality~\eqref{eq:lower-bound} holds.

\medskip
\noindent\textbf{Case 2.} $|H_G^*|\geq f+t$ and $(\bigcap_{v\in H_G^*}N_G(v))\setminus D=\varnothing$.

By Lemma~\ref{lem4}, there is a set $X\subseteq V(G)$ with $|X|\leq(f+t-2)^2$ such that every vertex outside $X$ has at least $t+1$ neighbors in $X$. Consequently, $\sum_{v\in X}d_G(v)\geq e_G(X,V(G)\setminus X)\geq(t+1)(n-|X|)$. By Lemma~\ref{lem3}, $\sum_{v\in V(G)\setminus X}d_G(v)\geq(t+2)(n-|X|)-|X|$. Thus, $e(G)\geq(t+\frac32)n-(t+2)(f+t-2)^2$. 
Note that  $n\geq 2(f+t-2)^3$. Thus, inequality \eqref{eq:lower-bound} holds.

\medskip
\noindent\textbf{Case 3.} $|H_G^*|\geq f+t$ and $(\bigcap_{v\in H_G^*}N_G(v))\setminus D\neq\varnothing$.

Choose $w$ from $(\bigcap_{v\in H_G^*}N_G(v))\setminus D$. Since $w$ is adjacent to every vertex of $H_G^*$, we have $w\in H_G$. By Lemma~\ref{lem6}, choose $L_G^*\subseteq N_G(L_G)$ such that $w\notin L_G^*$ and $|L_G^*|=|L_G|$. For every $z\in L_G^*$, there is a vertex $x\in L_G$ such that $w,z\in N_G(x)$. By Lemma~\ref{lem5}, $wz\in E(G)$. Thus $L_G^*\cup H_G^*\subseteq N_G(w)$. Moreover, Lemma~\ref{lem5} implies $L_G^*\subseteq H_G$, so $L_G^*\cap H_G^*=\varnothing$. Hence $d_G(w)\geq|H_G^*|+|L_G^*|$.  Note that $|L_G|=|L_G^*|$ and $|H_G^*|=|L_G|+|M_G|$. Thus $d_G(w)\geq 2|L_G|+|M_G|$. As $w\notin D$, there is a vertex $w'$ such that $ww'\in E(\overline G)$. Recall $H_G^*\subseteq N_G(w)$. Then $w'\in H_G$. 
We next establish three lower bounds on $e(G)$ and combine
them to obtain the desired inequality \eqref{eq:lower-bound}.

\begin{claim}\label{claim:estimatef1}
   $e(G)\ge \frac{1}{2} \big((t+3)|H_G|+(t+4)|M_G|+(t+4)|L_G|-(t+1)d_G(w')-2t-5\big)$.
\end{claim}
\begin{proof}
By Lemma~\ref{lem1},
$d_G(v,w')\le2$ for each $v\in V(G)\setminus\{w'\}$. Hence \(\sum_{v\in N_G(w')}d_G(v)\ge e_G(N_G(w'), V(G) \setminus N_G[w']) + e_G(N_G(w'),w')\geq n-1.\)
Let $R=V(G)\setminus (N_G[w']\cup\{w\})$. By Lemma~\ref{lem3},   \(\sum_{v\in R}d_G(v)\ge (t+2)(n-d_G(w')-2)-|L_G|.\) Recall $d_G(w)\geq 2|L_G|+|M_G|$.  Therefore, 
\[
\begin{aligned}
2e(G)
&= \sum_{v\in R}d_G(v)+d_G(w)+d_G(w')+\sum_{v\in N_G(w')}d_G(v)\\
&\ge (t+2)(n-d_G(w')-2)-|L_G|+|L_G^*|+|H_G^*|+d_G(w')+n-1\\
&=(t+3)|H_G|+(t+4)|M_G|+(t+4)|L_G|-(t+1)d_G(w')-2t-5. \quad \quad \quad \quad \quad \quad \quad\qedhere
\end{aligned}
\] 
\end{proof}
 
\begin{claim}\label{claim:estimatef2}
 $e(G)\ge \frac{1}{2} \big((f+t-2)|H_G|+(t+3)|M_G|+(t+4)|L_G|+d_G(w')-2f-2t+4\big)$.   
\end{claim}

 \begin{proof}
Since $w\notin L_G^*$, $L_G^*\subseteq N_G(w)$ and
$ww'\notin E(G)$, we have
$L_G^*\subseteq H_G\setminus\{w,w'\}$.
By Lemma~\ref{lem5}, every vertex of $L_G^*$ has degree at least
$f+t-1$. Hence $\sum_{v\in H_G\setminus\{w,w'\}}d_G(v)\geq(f+t-2)(|H_G|-2)+|L_G|$. Note that $\sum_{v\in H_G^*}d_G(v)\geq(t+2)|M_G|+(t+1)|L_G|$ and $d_G(w)\geq|M_G|+2|L_G|$.  Thus
\[
\begin{aligned}
2e(G)
&=\sum_{v\in H_G\setminus \{w,w'\}}d_G(v)+\sum_{v\in H_G^*}d_G(v)+d_G(w)+d_G(w')\\
&\ge (f+t-2)(|H_G|-2)+|L_G|+(t+2)|M_G|+(t+1)|L_G|+|M_G|+2|L_G|+d_G(w')\\
&=(f+t-2)|H_G|+(t+3)|M_G|+(t+4)|L_G|+d_G(w')-2f-2t+4.\quad \quad \quad\quad \quad \quad\qedhere
\end{aligned}
\]
 \end{proof}

 \begin{claim}
     $e(G)\ge |H_G|+(t+\frac{3}{2})|M_G|+(t+2)|L_G|-1-(f+t-3)(f+t-2)$. 
 \end{claim}
 
 \begin{proof}
For each $v\in V(G)$, define $wt(v)$ as follows. For each edge
$uv\in E(G)$, if $u,v\in H_G$, assign weight $0$ to both $u$ and $v$;
if $u\in H_G$ and $v\in H_G^*$, assign weight $1$ to $v$; and if
$u,v\in M_G$, assign weights $\frac{1}{2}$ and $\frac{1}{2}$ to $u$ and $v$, respectively.
Then $wt(v)$ is the total weight
assigned to $v$. By Lemma~\ref{lem5}, each neighbor of a vertex in $L_G$ has degree at least $f+t-1$. 
Hence, vertices in $L_G$ are adjacent only to vertices in $H_G$. 
Therefore, \(e(G)=e_G(L_G,H_G)+e_G(M_G,H_G)+e(G[M_G])+e(G[H_G])=\sum_{v\in V(G)}wt(v)+e(G[H_G]).\)

\begin{subclaim}\label{subclaim:high-neighbors}
    For each $v\in M_G$, the set $N_G(v)\cap H_G$ contains at least $t$ pairwise adjacent vertices.
\end{subclaim}
 \begin{proof}
Choose $v\in M_G$. Since $|H_G^*|\geq f+t$ and $d_G(v)\leq f+t-3$, there is a vertex $u\in H_G^*\setminus N_G[v]$. Since $G$ is $(K_t\vee \mathcal F)$-semi-saturated, $G+uv$ contains a copy of $K_t\vee F$ using $uv$ for some $F\in\mathcal F$. Every vertex of $H_G^*$ has degree at most $f+t-2$ in $G+uv$, so no vertex of $H_G^*$ can be a center of this copy. Thus, all $t$ centers belong to $N_G(v)\cap H_G$ and are pairwise adjacent.   \end{proof}

For each $v\in M_G$, since $d_G(v)\geq t+2$ and $|N_G(v)\cap H_G|\geq t$, we have $wt(v)=\frac12(d_G(v)+|N_G(v)\cap H_G|)\geq t+1$. Let $M_G^1$ denote the set of vertices in $M_G$ with weight $t+1$, and let $M_G^2$ denote $M_G\setminus M_G^1$.

\begin{subclaim}\label{subclaim:M1}
If $v\in M_G$ has at least $j$ neighbors in $M_G^1$,
where $j\geq1$ is an integer, then $vw'\in E(G)$ and
$wt(v)\geq\max\{t+2,t+1+j/2\}$.
In particular, $M_G^1$ is independent.
\end{subclaim}
\begin{proof}
Choose $u\in M_G^1\cap N_G(v)$.  By Claim~\ref{subclaim:high-neighbors},
$N_G(u)$ consists of $t$ pairwise adjacent vertices in $H_G$ and two vertices in
$M_G$, one of which is $v$. Write $U=N_G(u)\cap H_G$. Since $w\in U$  and $ww'\notin E(G)$, we have 
$w'\notin N_G(u)$.
As $G$ is $(K_t\vee\mathcal F)$-semi-saturated, $G+uw'$ contains a
copy $T$ of $K_t\vee F$ using $uw'$ for some $F\in \mathcal F$. Since each vertex of $H_G^*$ has degree at most $f+t-2$ in $G+uw'$, $u$ belongs to the $F$ part of $T$ and all centers of $T$ belong to $U\cup\{w'\}$.
Since $w'\in V(T)$ and $ww'\notin E(G)$, $w$ is not a center of $T$. Note that $|\{w'\}\cup U|=t+1$. Then $w'$ is a center and $w\notin V(T)$.  Since $u$ has at least two neighbors in the $F$ part of $T$
and $w\notin V(T)$, both of its neighbors in $M_G$
must belong to $T$. 
Hence $v\in V(T)$, 
$\{w'\}\cup U\subseteq N_G(v)$ and $wt(v)=|N_G(v)\cap H_G|+\frac12|N_G(v)\cap M_G|\geq t+1+\frac{j}{2}$. For $j\geq2$, this gives the required bound. Now assume that $j=1$. Suppose, to the contrary, that $wt(v)<t+2$. Then $v$ has exactly $t+1$ neighbors in $H_G$ and only one neighbor $u$ in $M_G$. Its neighbors in $H_G$ would be precisely the centers of $T$ and $w$. Since $w\notin V(T)$, the vertex $v$ would have only one neighbor $u$ in the $F$ part of $T$, contradicting $\delta(F)\geq2$.
Therefore, $wt(v)\geq t+2$. This also shows that every vertex with a neighbor
in $M_G^1$ lies outside $M_G^1$,
so $M_G^1$ is independent.  
\end{proof}

\begin{subclaim}
 We have $\sum_{v\in M_G}wt(v)\geq(t+\frac32)|M_G|$.
\end{subclaim}
\begin{proof}
Transfer $1/4$ of weight along each edge between $M_G^1$ and $M_G^2$, from its endpoint in $M_G^2$ to its endpoint in $M_G^1$. By Claim~\ref{subclaim:high-neighbors} and Claim~\ref{subclaim:M1}, every vertex of $M_G^1$ has exactly two neighbors in $M_G$, both of which belong to $M_G^2$. Its final weight is $t+\frac32$. For $v\in M_G^2$, let $j$ denote its number of neighbors in $M_G^1$. If $j=0$, then its weight is unchanged and is at least $t+\frac{3}{2}$.
If $j=1$, its final weight is at least $t+2-\frac14>t+\frac32$ by Claim 3.2. If $j\geq2$, its final weight is at least $t+1+\frac j2-\frac j4\geq t+\frac32$. Therefore, $\sum_{v\in M_G}wt(v)\geq(t+\frac32)|M_G|$.
\end{proof}

\begin{subclaim}
$e(G[H_G])\geq|H_G|+|L_G|-1-(f+t-3)(f+t-2)$.
\end{subclaim}
\begin{proof}
Let $G_1,G_2,\ldots,G_r$ denote the connected components of $G[H_G]$, and choose $x\in H_G^*$. By Lemma~\ref{lem1}, each component contains either a neighbor of $x$ or a neighbor of some vertex in $N_G(x)\cap H_G^*$. Since every vertex of $H_G^*$
has degree at most $t+f-3$, we have
$r\leq(t+f-3)+(t+f-3)^2=(t+f-3)(t+f-2)$. Since every connected component $G_i$ has at least $|V(G_i)|-1$ edges, we obtain $e(G[H_G])\geq|H_G|-r\geq|H_G|-(f+t-3)(f+t-2)$. If $|L_G|\leq1$, then
$e(G[H_G])
\geq|H_G|+|L_G|-1-(t+f-3)(t+f-2)$.
 Assume that $|L_G|\geq2$. By Lemma~\ref{lem5}, for each $y\in L_G$, $w$ is adjacent to every vertex of $N_G(y)\setminus\{w\}$. Since $ww'\notin E(G)$, we have $w'\notin N_G(y)$. Since $G$ is $(K_t\vee\mathcal F)$-semi-saturated, $G+yw'$ contains a copy $T$ of $K_t\vee F$ using $yw'$ for some $F\in\mathcal F$. The vertex $y$ has degree $t+2$ in $G+yw'$, so it belongs to the $F$ part of $T$ and $N_G(y)\subseteq V(T)$. In particular, both $w$ and $w'$ belong to the copy. Neither can be a center because they are nonadjacent. Thus, the center set is exactly $N_G(y)\setminus\{w\}$, and every vertex of this set is adjacent to $w'$. Recall $L_G^*\subseteq N_G(L_G)\setminus \{w\}$. 
Consequently, every vertex of $L_G^*$ is adjacent to both $w$ and $w'$, while $w'\notin N_G(L_G)$.
By Lemma~\ref{lem6}, $|L_G^*|=|L_G|$.
Thus, $G[H_G]$ contains a copy of $K_{2,|L_G|}$ with one
part $\{w,w'\}$.
We may assume that this copy is contained in $G_1$.
Since $G_1$ is connected, at least
$|V(G_1)|-(|L_G|+2)$ additional edges are outside this copy.
Thus,
$
e(G_1)\ge
2|L_G|+|V(G_1)|-(|L_G|+2)
=|V(G_1)|+|L_G|-2.
$
Since $e(G_i)\ge |V(G_i)|-1$ for $2\le i\le r$, we obtain
\begin{align*}
e(G[H_G])
\ge |V(G_1)|+|L_G|-2
   +\sum_{i=2}^r\bigl(|V(G_i)|-1\bigr)
\ge |H_G|+|L_G|-1-(f+t-3)(f+t-2).
\end{align*} 
\end{proof}

By the definition of $wt(v)$, we have  $\sum_{v\in L_G} wt(v)=\sum_{v\in L_G}d_G(v)=(t+1)|L_G|$.  
Together with Claims~3.3 and  3.4, we have

\[\begin{aligned}
e(G)=\sum_{v\in V(G)}wt(v)+e(G[H_G])
\ge
\left(t+\frac32\right)|M_G|
+|H_G|+(t+2)|L_G|-1-(f+t-3)(f+t-2).
\end{aligned}
\]
\end{proof}

Put $C=(t+1)(t+f-2)+t^2+t-2$. Taking a convex combination of the bounds in Claims~1--3
with coefficients
\[
\frac{2t+1}{C}, \qquad
\frac{(t+1)(2t+1)}{C}, \qquad
1-\frac{(t+2)(2t+1)}{C},\qquad
\]
respectively, which are positive since $f\ge8$, and using
$|H_G|+|M_G|+|L_G|=n$, we obtain
\[
e(G)>
\left(t+\frac32-
\frac{(2t+1)(t^2+2t-1)}
{2C}\right)n
-(t+f-2)^2.
\]
Moreover,
\[
\frac{t^2+2t-1}{t+f-2}
-\frac{(2t+1)(t^2+2t-1)}
{2C}
\ge\frac1{2(t+f-2)}.
\]
Together with $n\ge2(t+f-2)^3$, this gives
\[
\begin{aligned}
e(G)
&>
\left(t+\frac32-\frac{t^2+2t-1}{t+f-2}\right)n
+\frac{n}{2(t+f-2)}-(t+f-2)^2\\
&\ge
\left(t+\frac32-\frac{t^2+2t-1}{t+f-2}\right)n,
\end{aligned}
\]
as required.
Thus, inequality~(1) holds in all three cases. This completes the proof. 
\end{proof}

 \subsection*{Acknowledgements}
        Junxue Zhang was supported by the National Natural Science Foundation of China (No. 12601655).

\end{sloppypar}
\end{document}